%% file: NeurIPS20_PhysicsAware.tex
\documentclass{article}

\PassOptionsToPackage{numbers, compress}{natbib}

\usepackage[preprint]{neurips_2020}

\usepackage[utf8]{inputenc} 
\usepackage[T1]{fontenc}    
\usepackage{hyperref}       
\usepackage{url}            
\usepackage{booktabs}       
\usepackage{amsfonts}       
\usepackage{nicefrac}       
\usepackage{microtype}      

\usepackage{xcolor}
\usepackage{graphicx}
\usepackage{epstopdf, epsfig}
\usepackage{amsmath}
\usepackage{subcaption}
\usepackage{bm}
\usepackage{upgreek} 

\usepackage{algorithm,algpseudocode}
\usepackage{algorithm}
\usepackage{algpseudocode}

\def\mydepository{https://github.com/rmojgani/PhysicsAwareAE}

\title{Physics-aware registration based auto-encoder for convection dominated 
PDEs}

\author{%
  Rambod Mojgani\\
	Department of Aerospace Engineering\\
	University of Illinois at Urbana-Champaign\\
	Urbana, IL 61801 \\
	\texttt{mojgani2@illinois.edu} \\
   \And
	Maciej Balajewicz\\
	Department of Aerospace Engineering\\
	University of Illinois at Urbana-Champaign\\
	Urbana, IL 61801 \\
	\texttt{mbalajew@illinois.edu} \\
}

\begin{document}
\input{symbols}
\maketitle

\begin{abstract}
We design a physics-aware auto-encoder to specifically reduce the
  dimensionality of solutions arising from convection-dominated nonlinear
  physical systems.  Although existing nonlinear manifold learning methods seem
  to be compelling tools to reduce the dimensionality of data characterized by
  a large Kolmogorov n-width, they typically lack a straightforward mapping
  from the latent space to the high-dimensional physical space.  Moreover, the
  realized latent variables are often hard to interpret.  Therefore, many of
  these methods are often dismissed in the reduced order modeling of dynamical
  systems governed by the partial differential equations (PDEs).  Accordingly,
  we propose an auto-encoder type nonlinear dimensionality reduction algorithm.
  The unsupervised learning problem trains a diffeomorphic spatio-temporal
  grid, that registers the output sequence of the PDEs on a non-uniform
  parameter/time-varying grid, such that the Kolmogorov n-width of the mapped
  data on the learned grid is minimized.  We demonstrate the efficacy and
  interpretability of our approach to separate convection/advection from
  diffusion/scaling on various manufactured and physical systems.
\end{abstract}

\section{Introduction}

\input{introduction}

\section{Preliminaries}

\input{preliminaries}

\section{Physics-aware registration}

\input{method}

\section{Experiments}

\input{experiments}

\section{Conclusion}
\input{conclusion}

\medskip

\small
\bibliographystyle{unsrtnat}
\bibliography{NeurIPS20_PhysicsAware}{}

\end{document}

%% file: symbols.tex
\newcommand{\doc}{paper}
\newcommand{\mynote}[1]{ {\color{red} Note: #1 } } 
\newcommand{\apriori}{\textit{a priori}}
\newcommand{\aposteriori}{\textit{a posteriori}}
\newcommand{\nwidth}{n-width}
\newcommand{\offline}{off-line}
\newcommand{\online}{on-line}
\newcommand{\multistart}{multi-start}
\newcommand{\ndimensional}{n-dimensional}

\newcommand{\Pen}{\ensuremath{ \textit{Pe} }}
\newcommand{\Rey}{\ensuremath{ \textit{Re} }}

\newcommand{\dfunc}[1]{\ensuremath{ \dfuncnull \left(#1\right) }}
\newcommand{\dfuncnull}{\ensuremath{ \bm{f} }}
\newcommand{\statepar}{\ensuremath{ \bm{w} }}
\newcommand{\stateparinit}{\ensuremath{ \sysn{\statepar}{0}}}

\newcommand{\timecounterd}{\ensuremath{ n }} 

\newcommand{\sys}[1]{\ensuremath{ \bm{#1} }}
\newcommand{\sysT}[2]{\ensuremath{ \sys{#1}^{T} }}

\newcommand{\sysnscalar}[2]{\ensuremath{ {#1}{ \left[#2\right] } }}

\newcommand{\sysn}[2]{\ensuremath{ \sys{#1}{ \left[#2\right] } }}
\newcommand{\sysnC}[2]{\ensuremath{ \sys{#1}{ \left(#2\right) } }} 
\newcommand{\invsysn}[2]{\ensuremath{ \sys{#1}^{-1}{\left[#2\right]} }}
\newcommand{\pinvsysn}[2]{\ensuremath{ \sys{#1}^{+}{\left(#2\right]} }}
\newcommand{\sysnSup}[3]{\ensuremath{ \sys{#1}^{#3}{\left[#2\right]} }}
\newcommand{\sysnSub}[3]{\ensuremath{ \sys{#1}_{#3}{\left[#2\right]} }}
\newcommand{\sysnT}[2]{\ensuremath{ \sysnSup{#1}{#2}{T} }}
\newcommand{\dsysn}[2]{\ensuremath{ \dot{\sys{#1}}{\left(#2\right)} }}
\newcommand{\dsysnC}[2]{\ensuremath{ \dot{\sys{#1}}{\left(#2\right)} }}
\newcommand{\ddsysnC}[2]{\ensuremath{ \ddot{\sys{#1}}{\left(#2\right)} }}

\newcommand{\minimize}{\ensuremath{ \text{minimize} }}
\newcommand{\jmax}{\ensuremath{ j_{max} }}

\newcommand{\imag}[1]{\ensuremath{ \text{Im}\left\{ #1 \right\} }}
\newcommand{\real}[1]{\ensuremath{ \text{Re}\left\{ #1 \right\} }}
\newcommand{\realsetnull}{\ensuremath{ \mathbb{R} } } 
\newcommand{\realsetO}[1]{\ensuremath{ \realsetnull^{#1} } } 
\newcommand{\realsetp}{\ensuremath{ \realsetO{+} }}
\newcommand{\realset}[2]{\ensuremath{ \realsetO{#1 \times #2} }} 
\newcommand{\realsetT}[2]{\ensuremath{ \realsetO{#1 \times #2} }}
\newcommand{\realsetTh}[3]{\ensuremath{ \realsetO{#1 \times #2 \times #3} }}%
\newcommand{\complexsetnull}{\ensuremath{ \mathbb{C} } } 
\newcommand{\Trace}[1]{\ensuremath{ \Tr\left(#1\right) } } 
\newcommand{\Norm}[1]{\ensuremath{  \left\| #1 \right\| }}
\newcommand{\NormF}[1]{\ensuremath{  \Norm{#1}_F }}
\newcommand{\NormT}[1]{\ensuremath{  \Norm{#1}_2 }}
\newcommand{\abs}[1]{\ensuremath{  \left| #1 \right| }}
\newcommand{\inv}[1]{\ensuremath{  \left(#1\right)^{-1} }}
\newcommand{\pinv}[1]{\ensuremath{  \left(#1\right)^{+} }}
\newcommand{\eig}{\ensuremath{ \lambda_i }}
\newcommand{\eigs}[1]{\ensuremath{ \lambda_i \left( \bm{ #1 }\right) }}
\newcommand{\eigsr}[1]{\ensuremath{ \lambda_i \left( { #1 }\right) }}
\newcommand{\eigsn}[3]{\ensuremath{ \lambda_{\text{#3}} \left( { #1 
}\left(#2\right)\right) }}

\newcommand{\rankA}{\ensuremath{ r }}

\newcommand{\svds}[1]{\ensuremath{ \sigma_i \left( \bm{ #1 }\right) }}
\newcommand{\diag}[1]{\ensuremath{ \text{diag} \left( #1 \right) }}

\newcommand{\bigzero}{\mbox{\normalfont\Large\bfseries 0}}
\newcommand{\zerosnull}{\ensuremath{ \bm{0} }}
\newcommand{\zeros}[1]{\ensuremath{\zerosnull_{#1 \times #1} }}
\newcommand{\eye}[1]{\ensuremath{ \bm{I}_{#1 \times #1} }}
\newcommand{\ones}[2]{\ensuremath{ \onesnull_{#1 \times #2} }}
\newcommand{\onesnull}{\ensuremath{ \bm{1} }}
\newcommand{\transpose}[1]{\ensuremath{ {#1}^{T} }}
\newcommand{\normL}[1]{\ensuremath{ \mathcal{L}_{#1}} }
\newcommand{\normH}[1]{\ensuremath{ \mathcal{H}_{#1}} }
\newcommand{\manifold}{\ensuremath{ \mathcal{M} }}
\newcommand{\manifoldEul}{\ensuremath{ \manifold 	}}
\newcommand{\manifoldALE}{\ensuremath{ \manifold_\mathcal{F} 	}}
\newcommand{\manifoldS}{\ensuremath{ \mathcal{S} }}
\newcommand{\mapEultoALE}[1]{\ensuremath{ \mapEultoALEnull \left(#1\right) }}
\newcommand{\mapALEtoEul}[1]{\ensuremath{ \mapALEtoEulnull \left(#1\right) }}
\newcommand{\mapEultoALEnull}{\ensuremath{ \mathcal{G}  	}}
\newcommand{\mapALEtoEulnull}{\ensuremath{ \mathcal{G}^{-1} }}

\newcommand{\basisgridU}{\ensuremath{ \bm{U}_\textit{x} }}
\newcommand{\basisgridV}{\ensuremath{ \bm{V}_\textit{x} }}
\newcommand{\basisgridv}{\ensuremath{ \bm{v}_\textit{x} }}
\newcommand{\basisgridvn}{\ensuremath{ {\basisgridv}_{\timecounterd} }} 
\newcommand{\grido}{\ensuremath{\bm{x}_0  }}
\newcommand{\sizeGridROM}{\ensuremath{ r }}

\newcommand{\coarsenx}{\ensuremath{ \bm{S} }}
\newcommand{\coarsent}{\ensuremath{ \bm{P} }}
\newcommand{\coarsenM}{\ensuremath{ \bm{W} }}

\newcommand{\nnzX}[1]{\ensuremath{ \nnz \left(#1\right) }}
\newcommand{\nnz}{\ensuremath{ {nnz} }}
\newcommand{\error}{\ensuremath{ \varepsilon }}
\newcommand{\errort}{\ensuremath{ \error \left(t\right) }}
\newcommand{\errorInf}{\ensuremath{ \error_{\infty}  }}  
\newcommand{\sizeROM}{\ensuremath{ k }}
\newcommand{\snapshot}{\ensuremath{ \bm{M} }}
\newcommand{\snapshotr}{\ensuremath{ \widetilde{\snapshot} }}

\newcommand{\reducem}[1]{\ensuremath{ \widetilde{\bm{#1}}_{n} }}
\newcommand{\reduce}[1]{\ensuremath{ {\bm{#1}_\textit{r} }}}
\newcommand{\reducet}[2]{\ensuremath{ \sysn{\reduce{#1}}{#2} }}
\newcommand{\reducetC}[2]{\ensuremath{ \sysnC{\reduce{#1}}{#2} }}
\newcommand{\reducetT}[2]{\ensuremath{ \reduce{\bm{#1}}^{T} \left( 
#2 \right) }}

\newcommand{\trialb}{\ensuremath{ \bm{\Phi} }}
\newcommand{\testb}{\ensuremath{ \bm{\Psi} }}

\newcommand{\stateparred}{\ensuremath{ \reduce{\statepar} }}

\newcommand{\residual}[1]{\ensuremath{ \bm{R} \left(#1\right) }}


\newcommand{\STMd}[2]{\ensuremath{ \sysn{\phi_{A}}{#1, #2} }}
\newcommand{\STMdr}[3]{\ensuremath{ \sysn{\phi_#3}{#1, #2} }}
\newcommand{\cutsig}{\ensuremath{ \sigma_0 }}
\newcommand{\sizecutsig}{\ensuremath{ l }}

\newcommand{\stab}[1]{\ensuremath{ \hat{#1} }}  
\newcommand{\stabmat}[1]{\ensuremath{ \hat{\bm{#1}} }}  
\newcommand{\deltat}{\ensuremath{ \Delta n }}  
\newcommand{\sigmaD}{\ensuremath{ \sigma^\Delta }}

\newcommand{\conseps}{\ensuremath{ \epsilon }}

\newcommand{\Bc}{\ensuremath{ \sys{B}_\textit{c}} }  
\newcommand{\uc}{\ensuremath{ \sys{u}_\textit{c}} }

\newcommand{\Sct}{\ensuremath{ \sysn{\Sigma_\textit{c}}{\timecounterd} }} 
\newcommand{\Sc}{\ensuremath{ \bm{\Sigma}_\textit{c}}} 

\newcommand{\minimizey}{\ensuremath{ \sum_{\timecounterd=0}^{N_t} \NormT{ 
\sysn{y}{\timecounterd}- 
\sysn{\reduce{\stab{y}}}{\timecounterd}}  }}

\newcommand{\minimizeynt}[1]{\ensuremath{ \sum_{\timecounterd=0}^{#1} \NormT{ 
			\sysn{y}{\timecounterd}- 
			\sysn{\reduce{\stab{y}}}{\timecounterd}}  }} 
		
\newcommand{\fft}[1]{\ensuremath{ \textit{fft}\left( #1 \right) }}

\newcommand{\Dxx}{\ensuremath{ \bm{D}_{xx} }}  
\newcommand{\order}[1]{\ensuremath{ \mathcal{O} \left(#1\right) }}  
\newcommand{\tnfunc}{\ensuremath{ \sysnscalar{t}{\timecounterd} = \timecounterd 
\Delta t }}  

\newcommand{\frnn}[1]{\ensuremath{ f_{RNN} \left(#1\right) }}
\newcommand{\Nh}{\ensuremath{ k }}

\newcommand{\fencnull}{\ensuremath{ f_{\theta} }}
\newcommand{\fenc}[1]{\ensuremath{ \fencnull \left(#1\right)}}
\newcommand{\fdecnull}{\ensuremath{ g_{\theta} }}
\newcommand{\fdec}[1]{\ensuremath{ \fdecnull \left(#1\right)}}

\def\springcoef{\ensuremath{k}}

\newcommand{\std}{\ensuremath{\textit{SD}}}
\newcommand{\variance}{\ensuremath{\std^2}}
\newcommand{\mean}{\ensuremath{\mu}}
\newcommand{\disNormal}[2]{\ensuremath{ \mathcal{N}\left(#1, #2\right) }} 
\newcommand{\disUniform}[2]{\ensuremath{ \mathcal{U}\left(#1, #2\right) }} 

%% file: introduction.tex
Many physical phenomena are modeled using partial differential equations
(PDEs).  High accuracy simulations of these large-scale nonlinear systems often
require extensive computational resources.  Although cheaper and faster
processors, as well as highly parallel architectures, have made many of such
large-scale computations viable, reduced order models of these systems remain
attractive, especially in many-query simulations and real-time control.

The goal of reduced order modeling is to leverage the vast amount of data 
generated from high accuracy simulations to learn a low-dimensional model 
that can accurately and efficiently approximate the underlying dynamical system.
This is especially a challenging task for convection dominated PDEs, where 
the Kolmogorov~\nwidth\ of the snapshots of the solution is relatively large, 
i.e. the solution cannot be effectively reduced on a linear subspace.
Such problems emerge frequently in a broad range of applications, from 
Navier-Stokes equations (fluid dynamics) to Schr\"odinger equation 
(quantum-mechanical systems)~\citep{Mendible_arxiv_2019}.
In the machine learning community, the recognition of similar challenge dates 
back to 1990s and attempts in classification of handwritten 
digits~\citep{Hinton_NIPS_1995}, where presence of simple transformations such 
as translations and rotations in the data-set is well known to dramatically 
deteriorate the accuracy of linear methods such as principle component analysis 
(PCA). 
Fundamentally, other linear manifolds (subspaces) suffer from similar 
drawbacks,  examples include 
proper orthogonal decomposition (POD),
multidimensional scaling (MDS)~\citep{Cox_handbook_2008}, 
factor analysis~\citep{Friedman_statlearn_2001} and 
independent component analysis (ICA)~\citep{Friedman_statlearn_2001}.
Therefore, the high dimensionality of the data on any of these linear manifolds 
has incentivized a slew of nonlinear manifold learning approaches, such as
Iso-map~\citep{Tenenbaum_NIPS_1997},
kernel PCA~\citep{Mika_NIPS_1998}, 
locally linear embedding (LLE)~\citep{Roweis_science_2000}, 
Laplacian eigenmaps (LEM)~\citep{Belkin_NeuralC_2003},
semi-definite embedding (SDE)~\citep{Weinberger_ACM_2004},
auto-encoders~\citep{Hinton_science_2006},
t-SNE~\citep{Maaten_JMLR_2008}, and
diffeomorphic dimensionality reduction~\citep{Walder_NIPS_2008}.

Although many of the aforementioned nonlinear methods provide the sought after 
low-dimensional manifold, only a few provide the mapping from the learned 
low-dimensional to the high-dimensional manifold, for a 
survey~\citep[see][]{Lee_JCP_2020}.
This is especially important in reduced order modeling of PDEs, since the 
models 
are to be evolved in the parameters space or time on the low-dimensional 
manifold, i.e. evolving of the latent variables, and subsequently the latent 
variables have to be mapped to the physical high-dimensional manifold. 
Auto-encoders (AE), 
specifically convolutional auto-encoders (CAEs)~\citep{Masci_ICANN_2011} and 
deep convolutional generative adversarial networks 
(DCGANs)~\citep{Radford_CoRR_2015}, are amongst the successful methods used in 
dimensionality reduction of PDEs~\citep{Lee_JCP_2020,Cheng_JCP_2020}.
However, linear manifolds such as 
proper orthogonal decomposition (POD) and 
dynamic mode decomposition (DMD) are still often extensively preferred to these 
nonlinear approaches, since they provide an interpretable framework for 
analysis of the system, as well as controlling of the reduced system.
POD reveals the coherent structures in fluid flows~\cite{Noack_2011, 
Holmes_2012}, and 
DMD obtains a finite-dimensional, matrix approximations of the Koopman 
operator, which opens the possibility of taking advantage of estimation and 
control theories developed for the linear systems~\citep{Kutz_DMD_2016}.
In a more recent effort, it is shown that deep AE architectures can be trained 
to transform nonlinear PDEs into linear PDEs, by learning the eigen-function of 
the Koopman operator~\citep{Gin_arxiv_2019}. 
In this approach, although the transformation is nonlinear, the latent 
variables lie on a linear subspace.
Finally, a similar approach that prioritizes the optimal reducibility by 
learning a nonlinear manifold leads to a low-dimensional latent space.
Therefore, by definition, such an approach results in a more efficient reduced 
order model.

In this \doc, we develop an AE to learn a manifold on which the reduced order 
models can be efficiently constructed. 
To this end, we pose an unsupervised learning problem, that learns a 
spatio-temporal grid on which the low-rank linear decomposition of the solution 
of the PDE is optimal. 
The method can be interpreted as learning a map that registers the output 
sequence of a convection-dominated PDE to a low-rank reconstruction of the 
solution.
The method is in spirit of registration based manifold learning 
approaches, e.g.~\cite{Walder_NIPS_2008,Taddei_SIAM_2020}.
Furthermore, we demonstrate the efficacy of our approach in reducing some 
inherently challenging problems.
Finally, we discuss how the realized grid can be interpreted as a low-rank 
representation of information traveling in the domain and how it separates 
convection/advection from diffusion/scaling on the example problems.

%% file: preliminaries.tex
In approximation theory, Kolmogorov \nwidth\ is used to measure how well 
the data -- the solution of the PDEs in the context of this \doc -- can be 
approximated on a linear manifold, i.e. subspace.
The connection between the Kolmogorov \nwidth~ and POD of the data is 
rigorously established~\citep{Djouadi_PACC_2010}, leading to a measure on the
accuracy and feasibility of a low-rank reduced order model on a subspace.

Constructing the bases from snapshots in the spirit of the POD method can 
be formulated as a low-rank matrix approximation problem as follows: 
For a given snapshot matrix $\snapshot \in \realsetT{N}{K}$,
find a low rank matrix $\snapshotr \in \realsetT{N}{K}$ that solves the 
minimization problem
\begin{equation}
\label{Eqn:LRA}
\begin{aligned}
&\underset{{\rm rank}(\snapshotr)= \sizeROM}{\minimize}
&& \displaystyle \NormF{ \snapshot - \snapshotr },
\end{aligned}
\end{equation}
where for a successful reduction, $\sizeROM \ll N$ and $\sizeROM \ll K$. 
A snapshot matrix is here a matrix whose columns contain the states of the 
system of interest. 
More specifically, each column corresponds to the state of the system for some 
particular value of the system parameters, time, or the boundary/initial 
conditions.
Hence, $ \snapshot = \begin{bmatrix} \bm{w}_{1} , \cdots , \bm{w}_{K}
\end{bmatrix}$, where $\bm{w}_{i} \in \realsetO{N}$ is the state at 
the $i^{th}$ parameter or time step.

In problem~\eqref{Eqn:LRA}, the rank constraint can be
taken care of by representing the unknown matrix as 
$\snapshotr=\sys{U} \sys{V}$, where $\sys{U} \in \realsetT{N}{\sizeROM}$ and
$\sys{V} \in\realsetT{\sizeROM}{K}$, so that problem~(\ref{Eqn:LRA}) becomes
\begin{equation}
\label{Eqn:LRA2}
\begin{aligned}
&\underset{\sys{U}, \sys{V} }{\minimize}
&& \displaystyle \NormF{ \snapshot -  \sys{U}\sys{V} }.
\end{aligned}
\end{equation}
It is well known that the solution of the above low-rank approximation 
problem
is given by the singular value decomposition (SVD) of $\bm{M}$. 
Specifically, 
$\bm{U}= \left[\bm{u}_1, \cdots, \bm{u}_{k} \right] \in 
\realset{N}{\sizeROM}$ 
and 
$\bm{V}= \bm{\Sigma} \left[\bm{v}_1, \cdots, \bm{v}_{\sizeROM} \right] 
\in \realset{k}{K}$,
where 
$\bm{M}=\bm{U}^* \bm{\Sigma}^* \bm{V}^{*{\rm T}}$,
$\bm{U}^{*} = \left[\bm{u}_1, \cdots, \bm{u}_{\sizeROM}, 
\bm{u}_{\sizeROM+1}, 
\cdots, \bm{u}_N \right]$,
$\transpose{\bm{V}^{*}} = \left[\bm{v}_1, \cdots, \bm{v}_{\sizeROM}, 
\bm{v}_{\sizeROM+1}, \cdots, \bm{v}_K \right]$, and
$\bm{\Sigma}^* = \diag{\sigma_1, \sigma_2, \cdots, \sigma_{\rankA}}$ 
is a diagonal rank-$r$ matrix of singular values, where $\sigma_1 \geq 
\sigma_2 
\geq \cdots \geq \sigma_{\rankA}$,
and  
$\bm{\Sigma} = \diag{\sigma_1, \sigma_2, \cdots, \sigma_{\sizeROM}}$.
This is called ``method of snapshots'' for computing the POD 
bases~\citep{Sirovich_QAM_1987}.
This decomposition has a very close connotation to factor 
analysis~\citep{Friedman_statlearn_2001} and can be reproduced by artificial 
neural networks with linear activations~\citep{Colin_NP_1197}.
Although the linearity of the learned manifold leads to inefficiencies in
convection-dominated PDEs with large Kolmogorov \nwidth, the existence of a 
closed form solution as well as the abundance of computationally efficient 
approaches, such as~\citep{Holmes_NIPS_2009}, has made the method predominately 
be utilized in the field.
Our goal is to extend the norm minimization problem of~\eqref{Eqn:LRA2} to an
interpretable and efficient nonlinear manifold learning problem.

%% file: method.tex
We generalize the linear manifold learning problem of~\eqref{Eqn:LRA}, as a 
nonlinear manifold learning as follows:
For a given data-set lying on the high-dimensional manifold, learn a manifold, 
and the corresponding mapping, on which the solution can be expressed as a 
low-rank linear decomposition.
The map from the high-dimensional physical manifold to the learned manifold is 
denoted by $\mapEultoALEnull$ and the inverse by $\mapALEtoEulnull$.
In the finite-dimension matrix form where 
$\mapEultoALEnull: \realsetT{N}{K} \rightarrow \realsetT{N}{K}$ and 
$\mapALEtoEulnull: \realsetT{N}{K} \rightarrow \realsetT{N}{K}$, 
the minimization problem has the following 
form:
\begin{equation}\label{Eqn:LaMIN1}
\begin{aligned}
& \underset{}{\minimize}
& &  
\NormF{  
	\bm{M} - \mapALEtoEul{\snapshotr} 
}.
\end{aligned}
\end{equation}
In this case, $\snapshotr=\sys{U} \sys{V}$ is the linear decomposition the data 
on the learned manifold, where $\sys{U} \in \realsetT{N}{\sizeROM_r}$ and 
$\sys{V} 
\in\realsetT{\sizeROM_r}{K}$.
In principle, the compression of the data on the learned manifold is lossless.
For the compression to be successful $\sizeROM_r \ll N$ and $\sizeROM_r \ll K$, 
and it outperforms~\eqref{Eqn:LRA} 
if and only if $\sizeROM_r < k$.

In this formulation, both the high-dimensional snapshots and its 
low-dimensional representation on the learned manifold (latent variables) are 
stored in a matrix. 
In contrast to most  machine vision tasks dealing with images, the 
data does not necessarily lie on a uniform Cartesian grid.
This is important since many of the PDEs are discretized on unstructured 
computational grids, while many of the traditional machine learning tools 
are developed for uniform Cartesian grids representing an image.
Often, the snapshot matrix of $\snapshot$ is defined on a constant grid (in 
physics: Eulerian framework) and $\snapshotr$, by construct, is associated to 
the snapshots of latent variables on a parameter/time-varying grid (in physics: 
arbitrary Lagrangian Eulerian framework). 
Therefore, \eqref{Eqn:LaMIN1} can be interpreted as a registration task, that 
minimizes the Kolmogorov \nwidth\ of the snapshots of the latent variables on 
the learned parameter/time-varying grid.

\subsection{Diffeomorphism and interpolation}
\label{sec:diffeomorphism}

We impose diffeomorphism as a condition on the mapping to and from the learned 
manifold:
By definition, a map, $\mapEultoALEnull$, is said to be diffeomorphic if 
$\mapEultoALEnull$ and $\mapALEtoEulnull$ are 
differentiable~\citep{Modersitzki_2009}.
Bijectivity (i.e. one to oneness) and smoothness guarantee diffeomorphism. 
Therefore, by enforcing the map to be diffeomorphic, we ensure existence and 
uniqueness of $\snapshotr$, given a $\snapshot$ and vice versa.
Bijectivity is achieved by ensuring that volume of all the cells remain 
strictly positive. 
A negative cell volume leads to the indeterminate derivative of the state 
parameter, which can be seen as a ``tear'' in an image.
Smoothness of the grid is maintained by penalizing abrupt changes of 
the grid volume, both in space and parameter/time.
Any of the commonly used interpolation schemes between the constant grid and 
the learned grid can be utilized in this setting.

\subsection{Low-rank registration}

So far, we only have tied the mapping to a time/parameter-varying grid, and have
suggested that the mapping between the constant and the learned grid is simply 
an interpolation between them. 
In this section, the construction of the grid in the parameter/time space is 
discussed.

There are two general approaches to formulate the grid deformation in a 
registration problem.
In the first class of approaches, the grid nodes are controlled as the solution 
of the a minimization problem. 
Diffeomorphism can then be achieved by enforcing  a constraint on the 
determinant of the deformation gradient: be strictly positive for all grid 
cells.
This approach leads to a high-dimensional optimization problem which its
nonlinearity and ill-posedness makes it computationally 
challenging~\cite{Mang_OE_2018}.
In the second class of approaches, the mapping is the solution of a transport 
equation, i.e. flow fields, as in diffeomorphic dimensionality reduction 
\citep{Walder_NIPS_2008}.
Interestingly, in some special cases, a similar transport equation arises where 
the frame of references is changed from the Eulerian to the Lagrangian 
viewpoint, i.e. by solving the hyperbolic PDEs on the corresponding 
characteristic lines.
This change of the reference is proven to be efficient in reduced order 
modeling  of convection dominated PDEs~\citep{Mojgani_2017, Lu_JCP_2020}. 
The existence of a low-rank near-optimal grid for many of the 
convection-dominated PDEs are demonstrated; However, the extension of this 
change of frame for any arbitrary systems is not straightforward.
We leverage this premise in a data-driven setting. 
In our proposed method the coordinates of the  grid at the $n^\textit{th}$ 
parameter or time-step is constructed as:
\begin{equation}\label{Eqn:grid}
\sysn{x}{\timecounterd} := \basisgridU \basisgridvn,
\end{equation}
where $\basisgridvn \in \realsetO{r}$ is the $n^{\textit{th}}$ column of 
$\basisgridV$ and
$\basisgridU \basisgridV$ is a rank-$\sizeGridROM$ matrix, 
$\basisgridU \in \realsetT{N}{\sizeGridROM}$ and $\basisgridV \in 
\realsetT{\sizeGridROM}{K}$, that represents the evolution of the 
parameter/time-varying grid, on which the low-dimensional latent variables lie.
The rank-$\sizeGridROM$ grid is then learned in~\eqref{Eqn:LaMIN1}.
The latent variables can be interpreted as the evolution of the state 
parameters on the low-rank approximation of how the information travels, i.e. 
characteristic lines of the hyperbolic PDEs.
This is one of the key elements of the proposed method, which greatly 
reduces the size of the optimization compared to existing registration-based 
methods, such as~\cite{Taddei_SIAM_2020}. More importantly, it results in 
unprecedented predictive capabilities, i.e. accuracy beyond the training range.

\subsection{Implementation}
\label{sec:implementation}

In this section, we summarize the elements of the proposed algorithm and make 
some clarifications on implementation of the method.
As illustrated in Fig.~\ref{fig:idea}, the procedure is designed to learn a 
low-rank grid, $\basisgridU \basisgridV$, on which the snapshot of the data, 
$\snapshotr = \mapEultoALE{\snapshot}$, is low-rank.
\input{method_schematic}
The final minimization problem has the following form:
\begin{equation}\label{Eqn:ResMin_ALE_smooth}
\begin{aligned}
& \underset{\basisgridU,
			\basisgridV}{\minimize}
& &  
\NormF{  
	\bm{M} - \mapALEtoEul{\bm{U}\bm{V}} 
} + 
\NormF{ 
	\bm{\Gamma}_1   \basisgridU 
} +
\NormF{ 
	\basisgridV   \ \transpose{\bm{\Gamma}_2} 
}, \\
& \text{subject to}
& & \sysn{v}{n}\geq v_{\text{min}}, \text{ for } n=1,\ldots,K, \\
&&& \sysn{x}{n}|_{\partial \Omega} = \bm{x}|_{\partial \Omega}, \text{ for }
n=1,\ldots,K,
\end{aligned}
\end{equation}
where
${\bm{\Gamma}}_1 \in \realsetT{N}{N}$
and
${\bm{\Gamma}}_2 \in \realsetT{K}{K}$
are Tikhonov matrices designed to promote grid smoothness.
Also
$ \sysn{v}{\timecounterd} \in \realsetO{N-1}$
is a vector of cell volumes of the parameter/time-varying grid at the 
$n^{\text{th}}$ parameter/time step,
$v_{\min}$ is the minimum admissible cell volume, and
$\sysn{x}{\timecounterd} |_{\partial \Omega}$
and
$\bm{x}|_{\partial \Omega}$
are boundary points of the learned grid and the constant grid, respectively. 
The appropriate Tikhonov matrices and the minimum cell volumes are 
hyper-parameters and problem dependent.

Moreover, for practical reasons, a weak constraint on the rank reduction is 
chosen.
While~\eqref{Eqn:LaMIN1} implies minimizing over the rank of $\snapshotr$, in 
many cases, the solution of the minimization for a preset size of the 
decomposition is preferred.
We assume $\snapshotr=\sys{U} \sys{V}$, where
$\sys{U} \in \realsetT{N}{\sizeROM_r}$, 
$\sys{V} \in\realsetT{\sizeROM_r}{K}$, and
$\sizeROM_r \ll N $ , $\sizeROM_r \ll K $.
Moreover, to reduce the size of the optimization problem, $\basisgridU$ and 
$\basisgridV$ are uniformly down-sampled, however, the objective is 
evaluated on the fine grid.

To interpolate the snapshots between the two sets of grid, we simply utilize a 
$p$-degree polynomial interpolation scheme.
This choice innately incorporates a sparsity pattern into the mapping.
While the latent space representation of a data-point on the learned grid using 
a nearest-neighbor interpolation only requires one data-point on the constant 
grid, a $p$-degree polynomial interpolation requires $p-1$ entries of the input 
vector.
This, in principle, leads to a great reduction in the size the optimization 
problem compared to the traditional neural networks, where there is no 
\apriori~assumption on the structure of the connectivity between the nodes.

Finally, the proposed mapping can be utilized as an auto-encoder layer in a 
neural network (Fig.~\ref{fig:architecture}). 
This additional layer, greatly improves accuracy and training cost of neural 
network based models.
The algorithm to train the proposed physics-aware registration based 
auto-encoder is summarized in Alg.~\ref{alg:alg}.
\input{schematic_autoencoder}
\input{Alg}

%% file: method_schematic.tex
\begin{figure}
	\centering
	\def\FileSchematic{data/schematic}
	\includegraphics[scale=1]{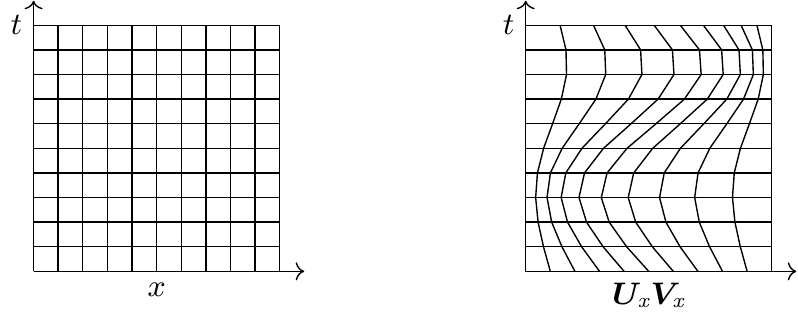}\\
	\hspace{-13pt}
	\includegraphics[scale=1]{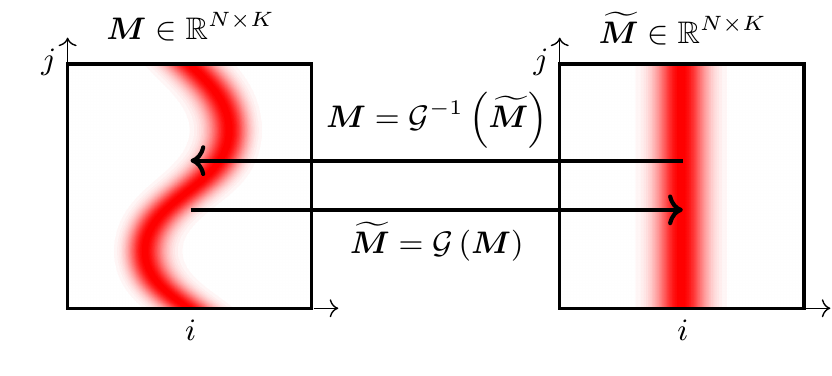}
	\caption{Illustration of proposed new nonlinear dimensionality reduction
	approach.}
	\label{fig:idea}
\end{figure} 
%

%% file: schematic_autoencoder.tex
\begin{figure}[!h]
	\centering
	\def\FileArch{data/RNN/}
	\begin{subfigure}[h]{0.39\textwidth}
		\centering
		\includegraphics[scale=0.25]{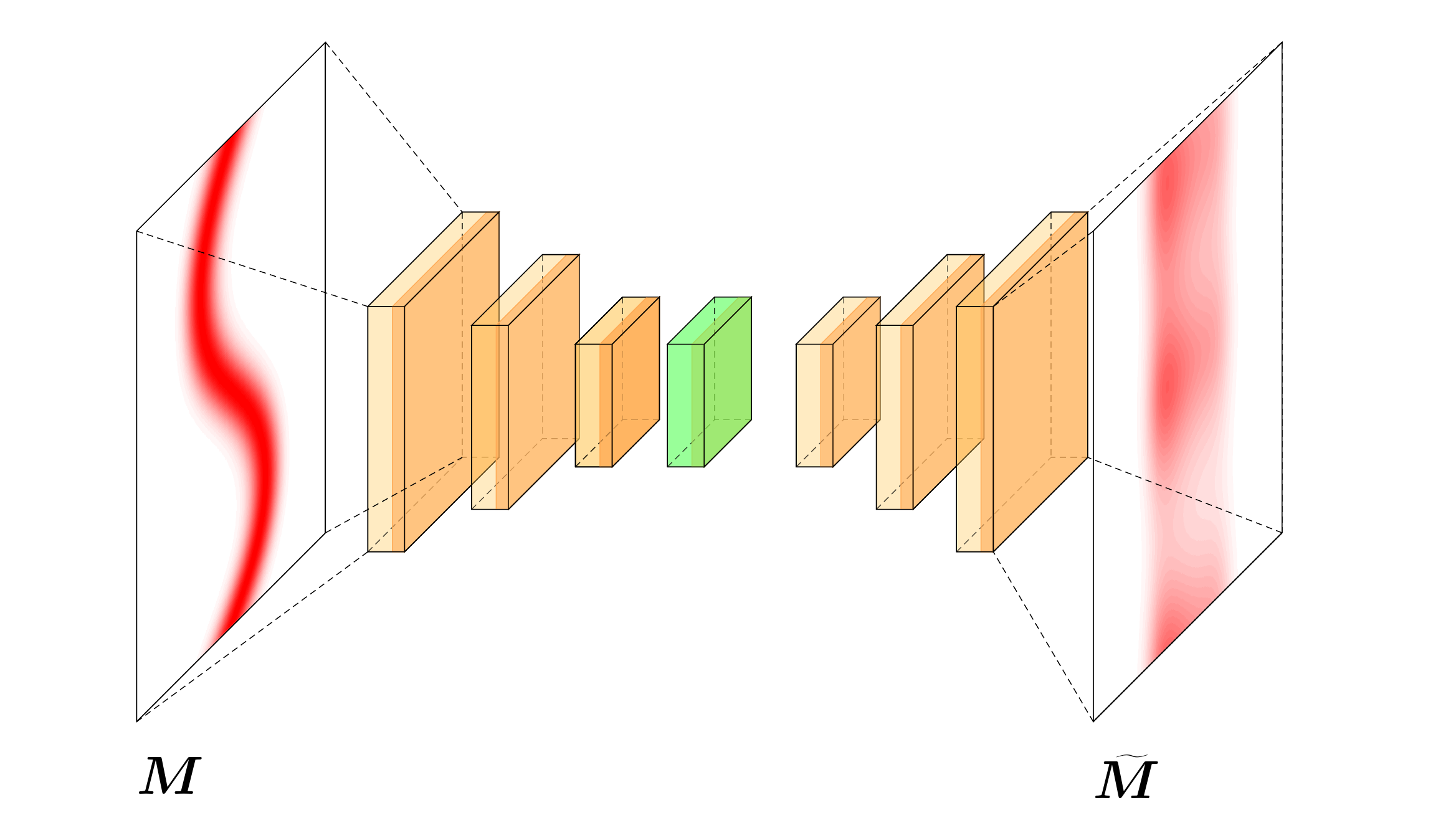}
		\caption{Traditional auto-encoder}
		\label{fig:arch_conventional}
	\end{subfigure}
	\begin{subfigure}[h]{0.59\textwidth}
		\centering
		\includegraphics[scale=0.25]{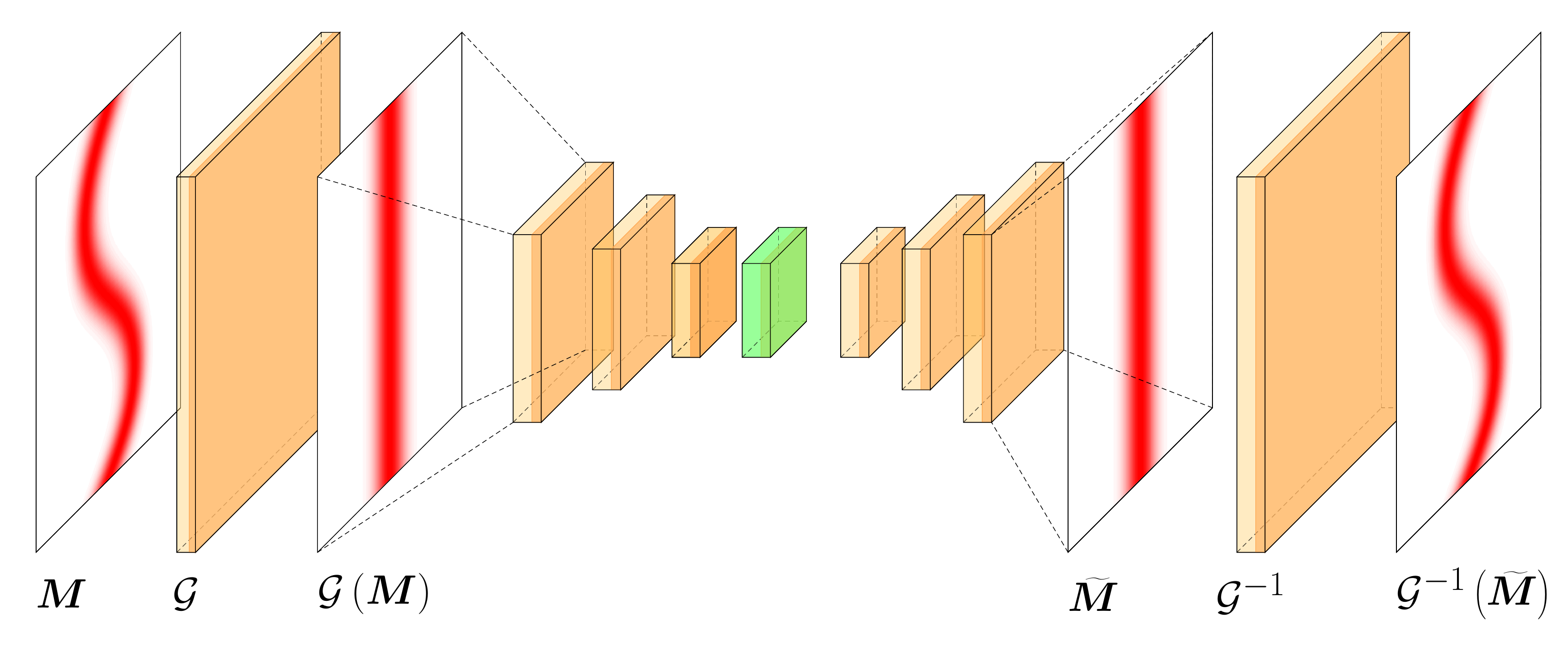}
		\caption{Physics-aware registration based auto-encoder}
		\label{fig:arch_physicsbased}
	\end{subfigure}
	\caption{
		Network architecture to map the solution of a dynamical system onto a
		low-dimensional space and map it back to the high-dimensional 
		space.
		The low-dimensional representation, the bottleneck, is depicted in 
		green.
	}
	\label{fig:architecture}
\end{figure} 

%% file: Alg.tex
\begin{algorithm}
		\caption{Training the physics-aware registration based auto-encoder}
		\label{alg:alg}
		\newcommand{\COMMENT}[2][.45\linewidth]{\small 
		\leavevmode\hfill\makebox[#1][l]{//~#2}}
	\hspace*{\algorithmicindent} \textbf{Input}:
	~
	Hyper-parameters: 
	\newline
	\hspace*{\algorithmicindent} \hspace*{\algorithmicindent}
	\hspace*{\algorithmicindent} \hspace*{\algorithmicindent}
	$\bm{\Gamma}_1$, $\bm{\Gamma}_2$,
	\newline
	\hspace*{\algorithmicindent} \hspace*{\algorithmicindent}
	\hspace*{\algorithmicindent} \hspace*{\algorithmicindent}
	Minimum admissible grid volume	($v_{\text{min}}$),
	\newline
	\hspace*{\algorithmicindent} \hspace*{\algorithmicindent}
	\hspace*{\algorithmicindent}
	Reduction parameters: 
	\newline
	\hspace*{\algorithmicindent} \hspace*{\algorithmicindent}
	\hspace*{\algorithmicindent} \hspace*{\algorithmicindent}
	Parameter/time-varying grid rank ($\sizeGridROM$),
	\newline
	\hspace*{\algorithmicindent} \hspace*{\algorithmicindent}
	\hspace*{\algorithmicindent} \hspace*{\algorithmicindent}
	Rank of the low-dimensional representation ($\sizeROM_r$),
		\newline
		\hspace*{\algorithmicindent} \hspace*{\algorithmicindent}
		\hspace*{\algorithmicindent}
		The snapshots matrix ($\snapshot \in \realset{N}{K}$),
		\newline
		\hspace*{\algorithmicindent} \hspace*{\algorithmicindent}
		\hspace*{\algorithmicindent}
		Maximum number of iterations ($j_{\text{max}}$),
		\newline
	\hspace*{\algorithmicindent} \textbf{Output}:
	\hspace{-4pt}
	Parameter/time-varying grid and its low-rank decomposition $\basisgridU 
	\basisgridV$, such that	
	$\sysn{x}{\timecounterd} := \basisgridU \basisgridvn$,
	and the corresponding maps, i.e. $\mapEultoALEnull$ and $\mapALEtoEulnull$
	\begin{algorithmic}[1]
		\State Initialize 
		$\basisgridU \in \realset{N}{\sizeGridROM}$ and 
		$\basisgridV\in \realset{\sizeGridROM}{K}$ with the SVD 
		decomposition of the constant grid plus a small random perturbation
		\State $j \leftarrow 0$
		\While{ $j \le j_\text{max}$ }
		\State $\snapshotr = \mapEultoALE{\snapshot}$ 
		\COMMENT{Interpolate the snapshots onto 
		$\basisgridU\basisgridV$}
		\State $\bm{U}\bm{V} \approx \snapshotr \text{~s.t.~} {\rm 
		rank}(\bm{U}\bm{V})= \sizeROM_r$ 
		\COMMENT{Approximate $\snapshotr$ using SVD as in~\eqref{Eqn:LRA2}}
		\State $\mapALEtoEul{\bm{U}\bm{V}}$ 
		\COMMENT{Interpolate $\bm{U}\bm{V}$ onto the constant grid}
		\State $\NormF{  
			\bm{M} - \mapALEtoEul{\bm{U}\bm{V}} 
		} +
		\NormF{ 
			\bm{\Gamma}_1   \basisgridU 
		} +
		\NormF{ 
			\basisgridV   \ \transpose{\bm{\Gamma}_2} 
		}$ 
		\COMMENT{Evaluate the objective}
		\State Update $\basisgridU$ and $\basisgridV$
		\COMMENT{Update the grid via the rule of the optimization}
		\State $j \leftarrow j+1$
		\EndWhile
	\end{algorithmic}
\end{algorithm}

%% file: experiments.tex
We provide several experiments to demonstrate the superior capabilities of the proposed 
physics-aware registration based auto-encoder. 
In this section, we solve for the low-rank grid and test the 
reducibility/compression of the input as well as its application in a neural 
network architecture approximate the PDEs.
We attempt to demonstrate both compression and interpretability of our 
results in these experiments.
For all these experiments, we resort to readily available optimization packages 
capable of solving optimization with nonlinear constraints; such as
\verb|minimize| in \verb|scipy| for our Python implementation, 
or interior-point method in \verb|fmincon| in our Matlab implementation.
Both of the implementations as well as the data and the results are made 
available at \url{\mydepository}.

\subsection{Manifold learning in rotated character ``A''}
\label{sec:rotated_A}

\begin{figure}[!h]
	\centering
	\def\rotatingA{data/RotatingA/n_mag_1_n_pos_2_edit}
	\begin{subfigure}[h]{0.99\textwidth}
		\centering
		\includegraphics[trim=200pt 610pt 140pt 
		75pt,clip,scale=0.35]{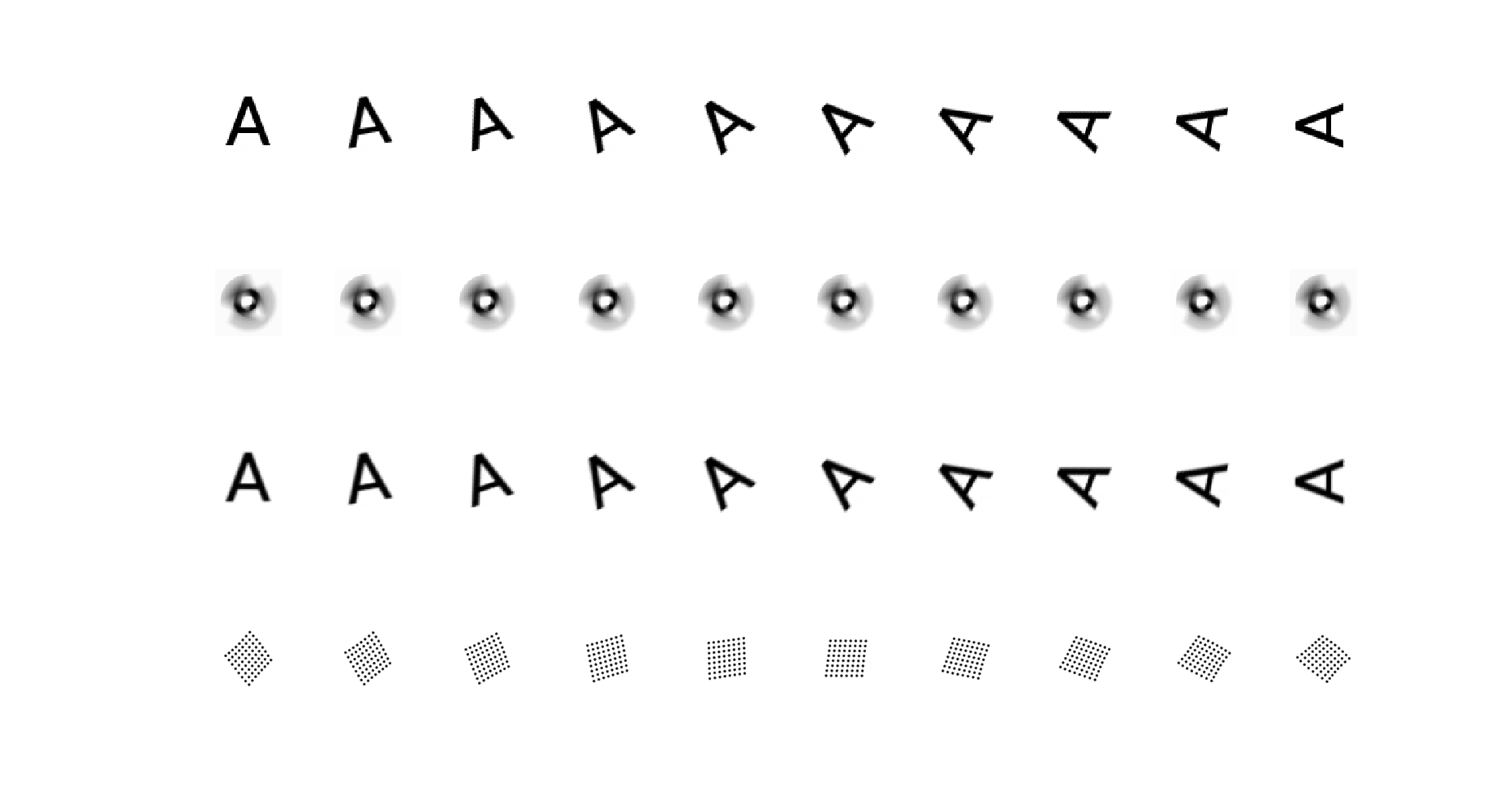}
		\caption{The snapshots.}
		\label{fig:roated_A_HFM}
	\end{subfigure}
	\begin{subfigure}[h]{0.99\textwidth}
		\centering
		\includegraphics[trim=200pt 445pt 140pt 250pt,clip,scale=0.35]
		{\rotatingA}
		\caption{The rank-$1$ reconstruction of the snapshots ($\sizeROM=1$) on 
		the POD subspace. }
		\label{fig:roated_A_EPOD}
	\end{subfigure}
	\begin{subfigure}[h]{0.99\textwidth}
		\centering
		\includegraphics[trim=200pt 110pt 140pt 570pt,clip,scale=0.35]
		{\rotatingA}
		\caption{The parameter-varying rank-$2$ grid corresponding to the 
			learned manifold.}
		\label{fig:roated_A_MG}
	\end{subfigure}
	\begin{subfigure}[h]{0.99\textwidth}
		\centering
		\includegraphics[trim=200pt 280pt 140pt 
		410pt,clip,scale=0.35]{\rotatingA}
		\caption{The rank-$1$ reconstruction of the snapshots ($\sizeROM_r=1$) 
		on the learned manifold. }
		\label{fig:roated_A_MGPOD}
	\end{subfigure}
	\caption{90 degrees rotation of character ``A''.}
	\label{fig:rotated_A}
\end{figure}

Consider a computer vision task of learning the nonlinear transformation, 
rotation, given a data-set comprised of a rotated character ``A''.  
The image of character ``A'' is stored in a $50 \times 50$ matrix and is 
rotated a total of 90 degrees with 3 degrees increments resulting in a snapshot 
matrix of dimension $2500 \times 31$.
A representative sample of the snapshots is shown in 
Fig.~\ref{fig:roated_A_HFM}, and a single POD mode reconstruction is 
illustrated in Fig.~\ref{fig:roated_A_EPOD}. 
In this problem, $\basisgridU$ is down-sampled to size of $7$, i.e. the total 
of $49$ control points. 
Moreover, $v_\text{min}=0$, $\bm{\Gamma}_1=100 \bm{D}_{x x}$ and 
$\bm{\Gamma}_2=\left(100/\pi\right) \bm{D}_{\theta \theta}$, where 
$\bm{D}_{x x}$ and
$\bm{D}_{\theta \theta}$ are the second derivative matrices in the spatial and 
parameter space, respectively.
The boundary point constraints are removed for this particular problem.
The optimization problem of~\eqref{Eqn:ResMin_ALE_smooth} approximates the 
rigid body rotation as in Fig.~\ref{fig:roated_A_MG}.
In Fig.~\ref{fig:roated_A_MGPOD}, the snapshots are approximated using a single 
basis ($\sizeROM_r=1$) on the learned manifold of $\sizeGridROM=2$.
The reconstruction delivered by the learned grid, using the proposed approach, 
is remarkably more accurate compared to the traditional POD approach.

\subsection{Manifold learning in two-dimensional fluid flows}
\label{sec:2DRiemann}

\begin{figure}[!h]
	\centering 
	\def\FileAddTDRie{data/2DRiemann/config03/t0d80/}
	\begin{subfigure}[!h]{0.32\textwidth}
		\includegraphics[scale=0.11]{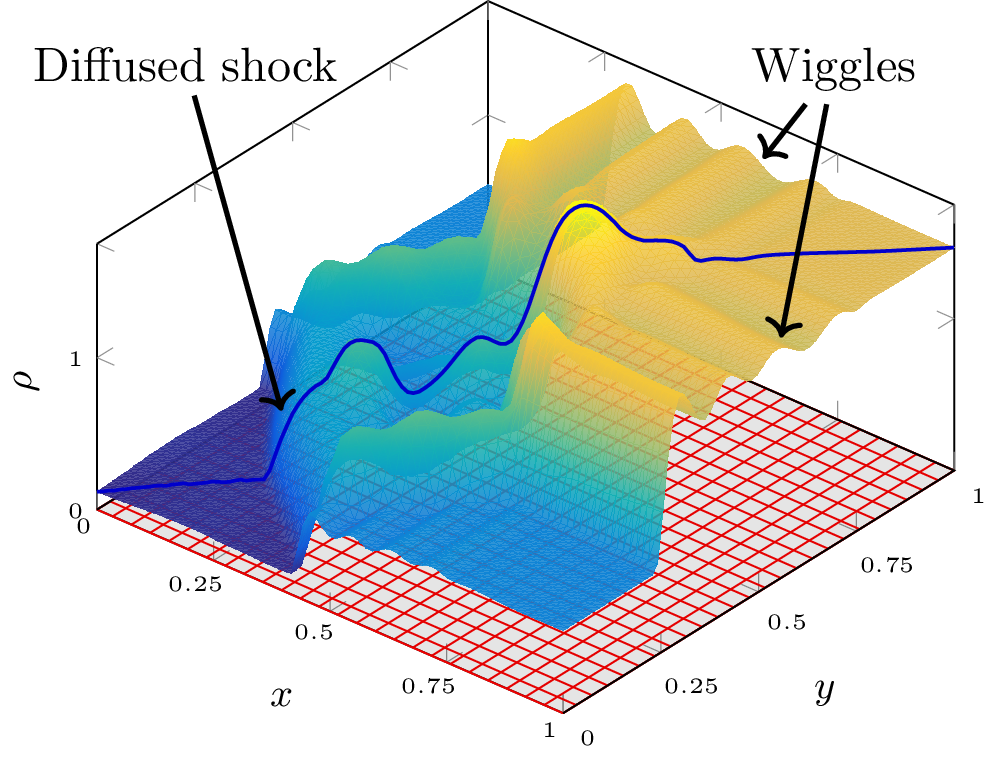}
		\caption{POD subspace($\sizeROM=8$)}
		\label{fig:2DRieROM_config03}
	\end{subfigure}
	\begin{subfigure}[!h]{0.32\textwidth}
		\includegraphics[scale=0.11]{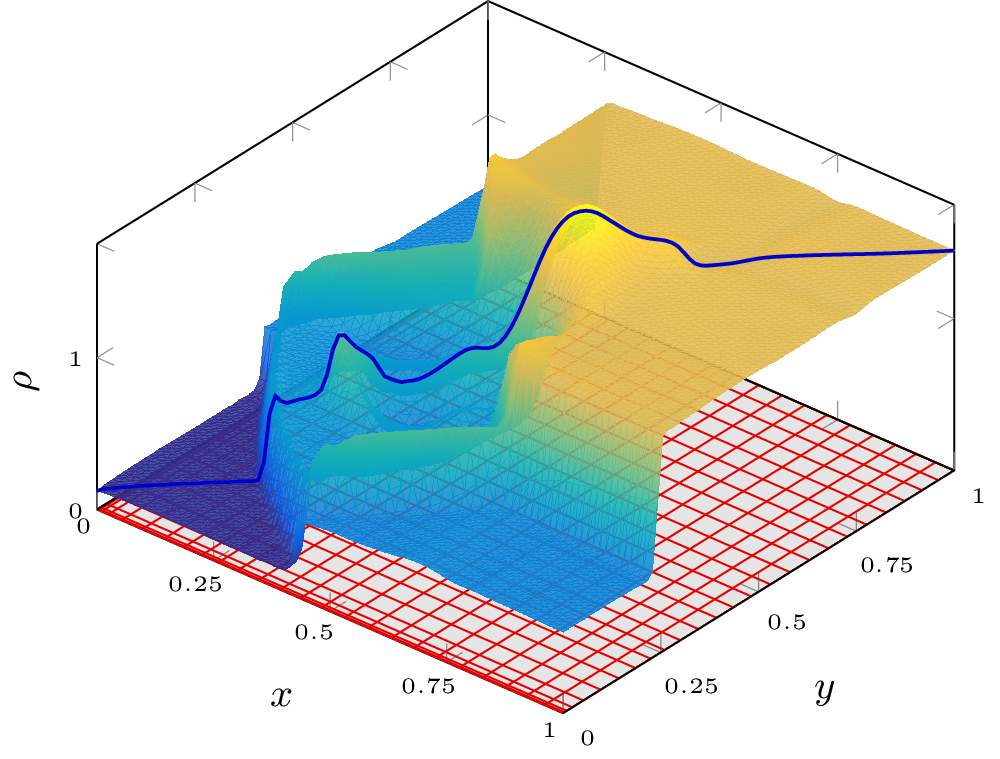}
		\caption{Learned manifold ($\sizeROM_r=8$)}
		\label{fig:2DRieMGROM_config03}
	\end{subfigure}
	\begin{subfigure}[!h]{0.32\textwidth}
		\includegraphics[scale=0.68]{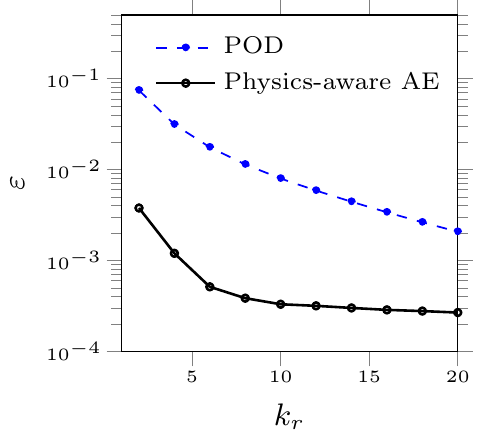}
		\caption{Error convergence}
		\label{fig:2DRie_config03_error}
	\end{subfigure}
	\centering
	\def\FileAddTDRie{data/2DRiemann/config12/t0d25/}
	\begin{subfigure}[h]{0.32\textwidth}
		\includegraphics[scale=0.11]{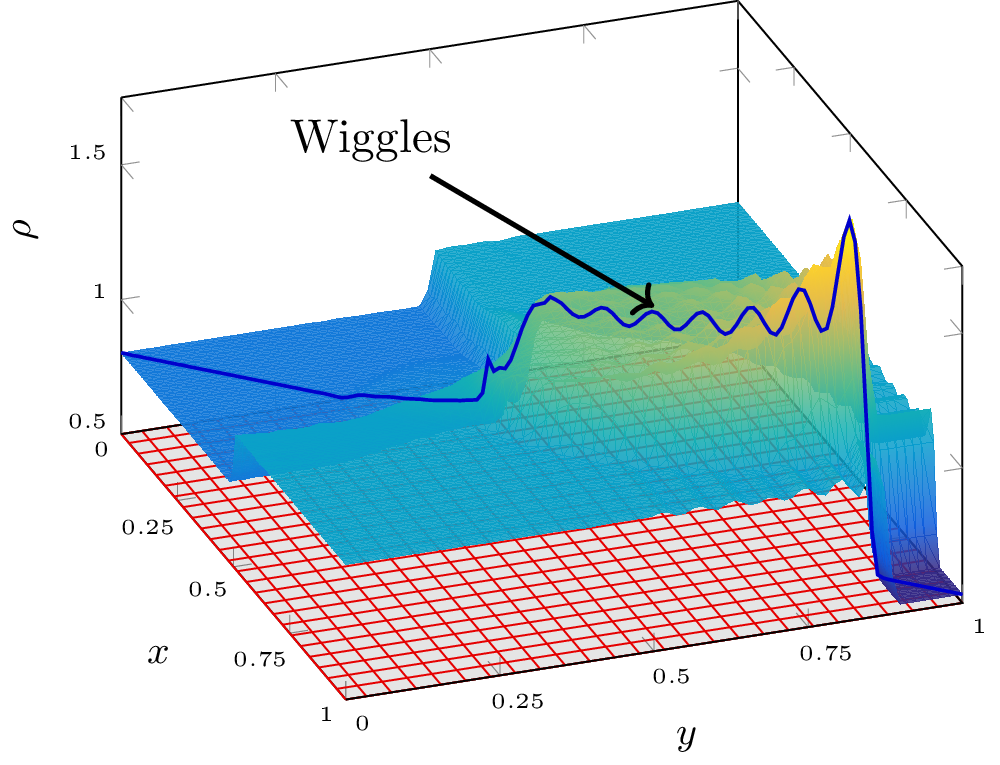}
		\caption{POD subspace ($\sizeROM=8$)}
		\label{fig:2DRieROM_config12}
	\end{subfigure}
	\begin{subfigure}[!h]{0.32\textwidth}
		\includegraphics[scale=0.11]{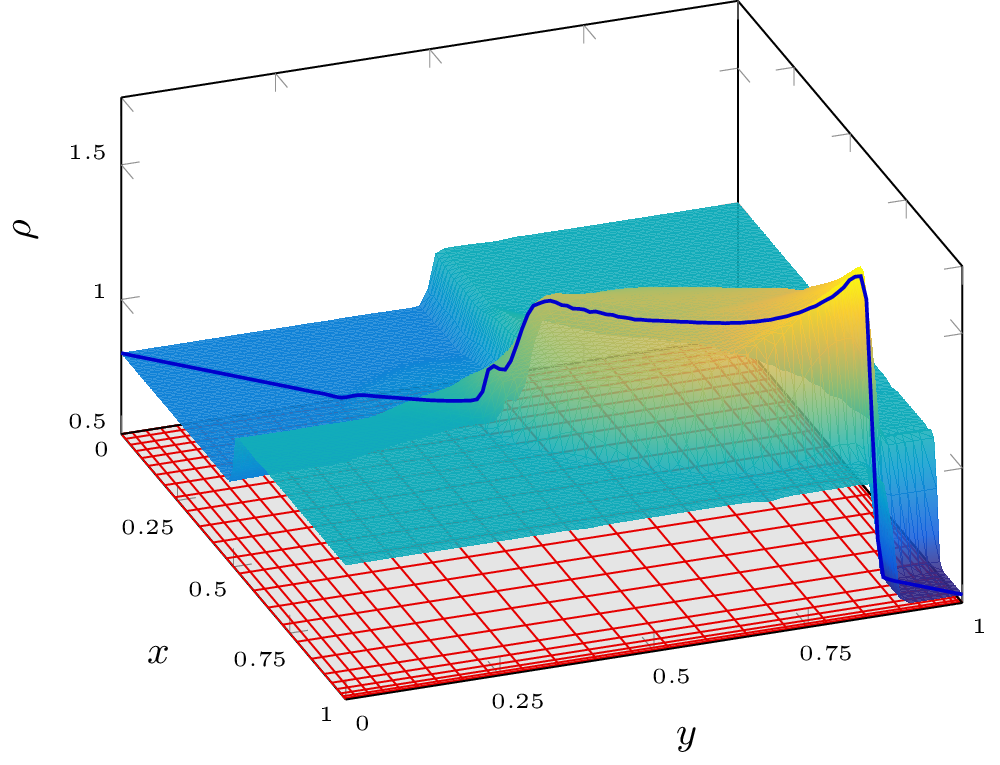}
		\caption{Learned manifold ($\sizeROM_r=8$)}
		\label{fig:2DRieMGROM_config12}
	\end{subfigure}
	\begin{subfigure}[!h]{0.32\textwidth}
		\includegraphics[scale=0.68]{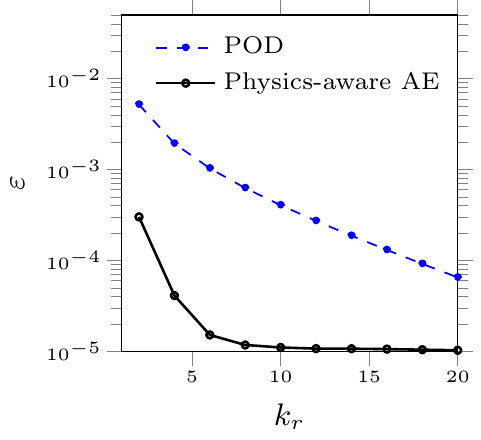}
		\caption{Error convergence}
		\label{fig:2DRie_config12_error}
	\end{subfigure}
	\caption{The density snapshots of the two-dimensional Riemann problems at 
	the last time step of the simulations and the computational grid 
	representative of the learned manifolds.
	Configuration 3: Fig.~\ref{fig:2DRieROM_config03} to 
	Fig.~\ref{fig:2DRie_config03_error};
	Configuration 12: Fig.~\ref{fig:2DRieROM_config12} to 
	Fig.~\ref{fig:2DRie_config12_error};}
	\label{fig:2DRie}
\end{figure}

In this section, the proposed manifold learning approach is applied to the 
two-dimensional Riemann problem of fluid dynamics,
$
\frac{\partial}{\partial t} \bm{q} + 
\frac{\partial}{\partial x} \bm{f}_x +
\frac{\partial}{\partial y} \bm{f}_y = 0,
$
where $\bm{q} = [\rho, \rho u, \rho v, \rho e]$,
$\bm{f}_x = [\rho u , \rho u ^ 2 + p, \rho u v , \rho u H]$,
$\bm{f}_y = [\rho v , \rho u v + p, \rho v^2 + p , \rho v H]$,
and $H = e + p/\rho$ ,
$p = \rho (\gamma -1 ) (e - 0.5 (u^2 + v^2) )$
in the domain $(x,y,t) \in [0,1] \times [0,1] \times [0,t_\text{max}]$, with 
initial conditions of configuration 3 and 12 as in~\cite{Lax1998}.
The snapshots of primitive variables are generated using a high-order 
artificial viscosity scheme coupled with a $4^{\textit{th}}$-order Runge-Kutta 
time discretization with $\Delta t = 5 \times 10^{-4}$ on a $150 \times 150$ 
grid.  Snapshots are collected at $\delta t=0.016$ and $\delta t=0.006$  
intervals for simulation range of $t \in [0,0.80]$ and $t \in [0,0.25]$
for configuration 3 and 12, respectively. 
The size of the snapshot matrices in both cases are $10000 \times 50$. 
A rank-$2$ time-varying grid ($\sizeGridROM=2$) is learned 
via~\eqref{Eqn:ResMin_ALE_smooth} setting $\sizeROM_r=4$, and 
$\gamma_1 = \gamma_2 = 0.05$,
where $\bm{\Gamma}_1 = \gamma_1 \bm{D}_{x x} = \gamma_1 \bm{D}_{y y}$ 
is the second derivative matrix in the $x$ and $y$ directions and 
$\bm{\Gamma}_2=\gamma_2  \bm{D}_{t t}$, where $\bm{D}_{t t}$ is the second 
derivative matrix in time.
Also, $v_\text{min}=\Delta x_\text{min}\Delta y_\text{min}$,
where $\Delta x_\text{min}=\Delta y_\text{min}=6.7 \times 10^{-4}$. 
The low-dimensional representation of density is compared on the constant grid 
and the learned manifold in Fig.~\ref{fig:2DRie}. 
On the learned manifold, traveling shocks are conserved and
free of non-physical oscillatory solutions, resulting in a 
significant error reduction.

\subsection{Physics based auto-encoder in an LSTM architecture}
\label{sec:RNN1D}

In this section, the proposed method is implemented as an auto-encoder layer 
wrapped around a traditional machine learning architecture to decrease the loss 
in compression phase and subsequently to improve the predictive capabilities of 
a recurrent neural network modeling the governing PDEs.
We employ a long short-term memory (LSTM) to approximate the dynamics of the 
PDE on the learned manifold, a proven architecture to approximate 
PDEs~\citep{Parish_CMAME_2020}. 
To learn a low-dimensional manifold, a dense neural network with 
linear activation is used. Finally, the proposed registration based 
auto-encoder is added to the architecture.

\begin{figure}[!h]
	\centering
	\def\FileRNN{data/RNN/burgers}
	\begin{subfigure}[h]{0.32\textwidth} 
		\includegraphics[trim={0 0 0 0},clip, scale=0.9]
		{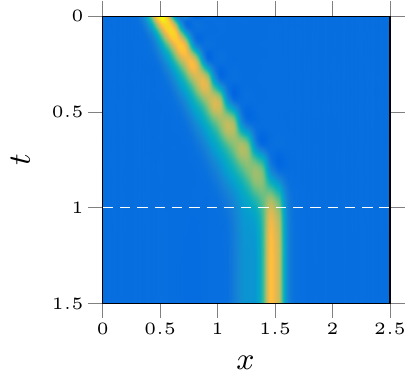}
		\caption{Neural Network AE ($N_h=10$)}
		\label{fig:Burgers_AERNN}
	\end{subfigure}
	\begin{subfigure}[h]{0.32\textwidth} 
		\includegraphics[trim={0 0 0 0},clip, scale=0.9]
		{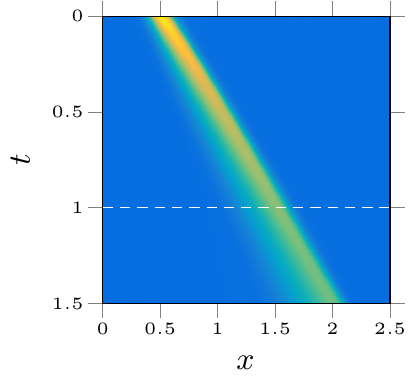}
		\caption{Physics-aware AE ($N_h=10$)}
		\label{fig:Burgers_PAAERNN}
	\end{subfigure}
	\begin{subfigure}[h]{0.32\textwidth} 
		\includegraphics[trim={0 0 0 0},clip, scale=0.78]
		{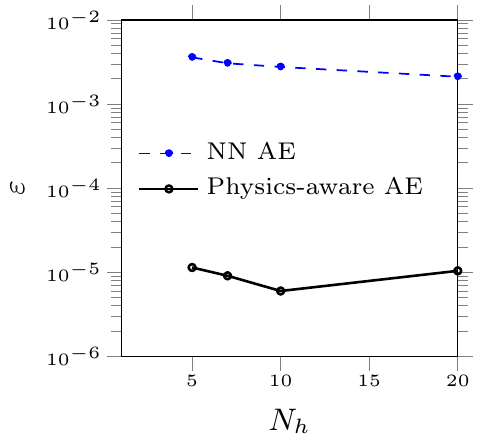}
		\caption{LSTM error}
		\label{fig:Burgers_RNN_error}
	\end{subfigure}
	\caption{LSTM approximating the Burgers' equation.} 
	\label{fig:Burgers_RNN}
\end{figure}

Consider the scalar, one-dimensional nonlinear convection-diffusion equation, 
known as viscous Burgers' equation,
$
\partial_t w(x,t)
+ w \partial_x w(x,t)
= \left(1/\Rey\right) \partial_{xx} w(x,t)
$
in the domain $(x,t) \in [x_a,x_b] \times [0,T]$, equipped with initial
conditions $w(x,0) = w_0(x)$, and Dirichlet boundary conditions at $x_a$, and
$x_b$, where $w(x,0)=0.8 + 0.5~e^{-(x-0.5)^2/0.1^2}$, $w(x_a,t)=w(x_b,t) = 0$, 
for $(x,t) \in [0,2.5] \times [0,1]$. 
An implicit second order time discretization is used with $\Delta 
t=8\times10^{-3}$ and space is uniformly discretized where $\Delta 
x=1\times10^{-2}$.
In the proposed architecture, the rank-$1$ time-varying grid 
($\sizeGridROM=1$), representing the physics-aware auto-encoder, is learned as 
in~\eqref{Eqn:ResMin_ALE_smooth} with $\sizeROM_r=4$.
In this problem, $v_\text{min}=10^{-3}$, $\bm{\Gamma}_1=\gamma_1 
\bm{D}_{xx}$ and $\bm{\Gamma}_2=\gamma_2 \bm{D}_{tt}$, where 
$\gamma_1=\gamma_2=1$, and $\bm{D}_{xx}$ and 
$\bm{D}_{t t}$ are second derivative matrices in space and time.
Subsequently, the LSTM is trained to approximate $\mapEultoALE{\snapshot}$.
In the traditional architecture, the densely connected auto-encoder and the 
LSTM are trained simultaneously (Fig.~\ref{fig:arch_conventional}).
In the proposed architecture (Fig.~\ref{fig:arch_physicsbased}), the manifold 
and the neural network are trained separately.
The neural network architectures are deployed in Keras~\citep{chollet2015keras}.
The models are trained using the snapshot matrix of 
$\snapshot\in\realset{N_x}{N_t}$.
In the proposed physics-aware AE, $\basisgridV$ is linearly extended. 
The contour plots of the solutions are compared in Fig.~\ref{fig:Burgers_AERNN} 
and Fig.~\ref{fig:Burgers_PAAERNN}.
As expected, the neural network AE cannot predict the convection underlying the 
physics of the problem outside the training range 
(Fig.~\ref{fig:Burgers_AERNN}); however, by levering the 
physics of the problem, the LSTM trained on the low-dimensional manifold 
realized by the registration-based AE leads to a solution much closer to the 
solution of the PDE (Fig.~\ref{fig:Burgers_PAAERNN}).
The mean squared difference of the solutions, $\error$, for different sizes of 
LSTM, $N_h$, are plotted in Fig.~\ref{fig:Burgers_RNN_error}.

\begin{figure}[!h]
	\centering
	\def\FileRNN{data/RNN/wave}
	\begin{subfigure}[h]{0.32\textwidth} 
		\includegraphics[trim={0 0 0 0},clip, scale=0.9]
		{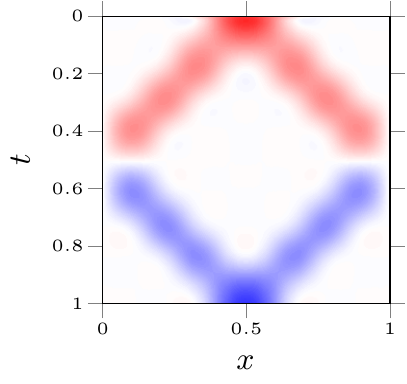}
		\caption{Neural Network AE ($N_h=5$)}
		\label{fig:wave_AERNN}
	\end{subfigure}
	\begin{subfigure}[h]{0.32\textwidth} 
		\includegraphics[trim={0 0 0 0},clip, scale=0.9]
		{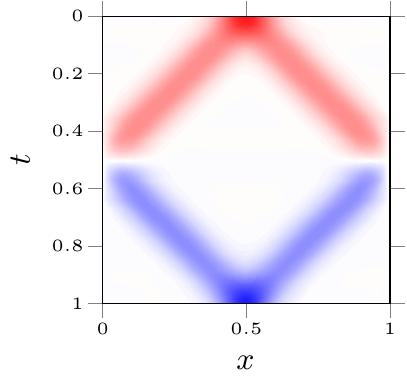}
		\caption{Physics-aware AE ($N_h=5$)}
		\label{fig:wave_PAAERNN}
	\end{subfigure}
	\begin{subfigure}[h]{0.32\textwidth} 
		\includegraphics[trim={0 0 0 0},clip, scale=0.78]
		{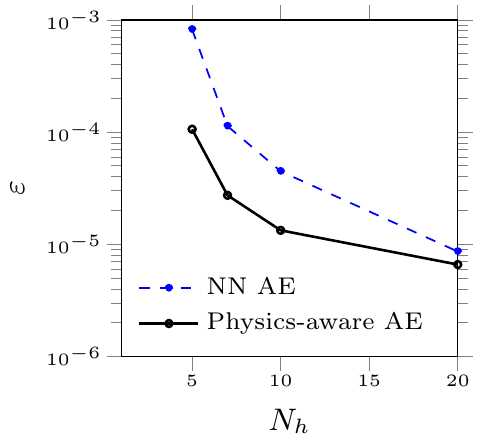}
		\caption{LSTM error}
		\label{fig:wave_RNN_error}
	\end{subfigure}
	\caption{LSTM approximating the wave equation.} 
	\label{fig:wave_RNN}
\end{figure}

Also, consider the one-dimensional wave equation,
$
\partial_{tt} w(x,t) - \partial_{xx} w(x,t) = 0,
$
in the domain $(x,t) \in [x_a,x_b] \times [0,T]$, equipped with initial
conditions $w(x,0) = w_0(x)$, and Dirichlet boundary conditions at $x_a$, and
$x_b$, where $w(x,0)= e^{-(x-0.5)^2/0.1^2}$, $w(x_a,t)=w(x_b,t) = 0$, 
for $(x,t) \in [0,1] \times [0,1]$. 
An implicit second-order time-discretization is used with $\Delta 
t=2.5\times10^{-3}$ and space is uniformly discretized where $\Delta 
x=10^{-2}$.
The architecture is set up similar to the Burgers' equation, with the following 
parameters:
the time-varying grid is of rank-$2$ ($\sizeGridROM=2$),
the reconstruction on the learned manifold is of rank-$2$ ($\sizeROM_r=2$), 
$v_\text{min}=10^{-3}$, $\gamma_1=\gamma_2=10$, and the size of the grid bases, 
in both space and time, are down-sampled to $15$ control points.
The solution and error are plotted in Fig.~\ref{fig:wave_RNN}, showing the 
performance of the proposed physics-aware AE compared to the traditional 
architecture.

%% file: conclusion.tex
We have proposed a variation to the diffeomorphic~\citep{Walder_NIPS_2008} or 
the registration based \citep{Taddei_SIAM_2020} manifold learning methods for 
dimensionality reduction of convection dominated PDEs, where the 
Kolmogorov~\nwidth\ of the snapshots is large.
The existence of a low-rank parameter/time-varying grid~\citep{Mojgani_2017, 
Lu_JCP_2020}, is leveraged to decrease the 
Kolmogorov~\nwidth\ of such problems in a data-driven settings.
In the proposed registration, the solution of the PDEs is mapped onto a 
non-uniform parameter/time-varying grid, such that the Kolmogorov \nwidth\ of 
the latent snapshots is minimized.
This is achieved using a low-rank linear decomposition of the 
snapshots on the learned grid, a relatively low-cost computation compared to 
training of any of the existing nonlinear manifold learning approaches.
We have demonstrated the application of the proposed physics-aware auto-encoder 
on a computer vision task, an image under nonlinear transformation, and 
solution of nonlinear convection dominated PDEs, two-dimensional Riemann 
problem and Burgers' equation.
Leveraging the low-rank structure of the grid has enabled us with an 
interpretable map that can be manipulated to fit the underlying convection in 
the problem.
We have incorporated the dimensionality reduction as an auto-encoder layer 
in training of recurrent neural networks reproducing PDEs.
Moreover, we have shown it outperforms a neural network-based auto-encoder in 
prediction of the system beyond the training range.
The proposed approach has the potential to be generalized to PDEs of higher 
spatial dimensions, as well as to greatly decrease the size of more 
traditional projection-based model order reduction techniques.
More importantly, it shows the potential to lower the training costs and 
enhance the predictive capabilities of models learned in system 
identification settings.